\documentclass[12pt]{article}
\usepackage{amssymb,amsthm,amsmath,latexsym}
\usepackage{url}
\usepackage{lscape} 
\newtheorem{thm}{Theorem}
\newtheorem{prop}[thm]{Proposition}

\theoremstyle{remark}

\newcommand{\FF}{\mathbb{F}}
\newcommand{\ZZ}{\mathbb{Z}}

\newcommand{\0}{\mathbf{0}}
\newcommand{\1}{\mathbf{1}}

\DeclareMathOperator{\wt}{wt}

\begin{document}

\title{On the minimum weights of binary linear complementary dual codes}

\author{
Makoto Araya\thanks{Department of Computer Science,
Shizuoka University,
Hamamatsu 432--8011, Japan.
email: {\tt araya@inf.shizuoka.ac.jp}}
and 
Masaaki Harada\thanks{
Research Center for Pure and Applied Mathematics,
Graduate School of Information Sciences,
Tohoku University, Sendai 980--8579, Japan.
email: {\tt mharada@tohoku.ac.jp}}
}

\maketitle

\noindent
{\bf Dedicated to Professor Masaaki Kitazume on His 60th Birthday}

\begin{abstract}
Linear complementary dual codes (or codes with complementary duals)
are codes whose intersections with their dual codes are trivial.
We study the largest minimum weight $d(n,k)$ among all binary linear 
complementary dual $[n,k]$ codes.
We determine $d(n,4)$ for 
$n \equiv 2,3,4,5,6,9,10,13 \pmod{15}$, 
and $d(n,5)$ for $n \equiv 3,4,5,7,11,19,20,22,26 \pmod{31}$.
Combined with known results, 
the values $d(n,k)$ are also determined 
for $n \le 24$.
\end{abstract}

%%%%%%%%%%%%%%%%%%%%%%%%%%%%%%%%%%
\section{Introduction}
Let $\FF_q$ denote the finite field of order $q$,
where $q$ is a prime power.
An $[n,k]$ {\em code} $C$ over $\FF_q$ 
is a $k$-dimensional vector subspace of $\FF_q^n$.
A code over $\FF_2$ is called {\em binary}.
The parameters $n$ and $k$
are called the {\em length} and {\em dimension} of $C$, respectively.
The {\em weight} $\wt(x)$ of a vector $x \in \FF_q^n$ is
the number of non-zero components of $x$.
A vector of $C$ is called a {\em codeword} of $C$.
The minimum non-zero weight of all codewords in $C$ is called
the {\em minimum weight} $d(C)$ of $C$. An $[n,k,d]$ code
is an $[n,k]$ code with minimum weight $d$.
Two $[n,k]$ codes $C$ and $C'$ over $\FF_q$ are 
{\em equivalent} %, denoted $C \cong C'$,
if there is an $n \times n$ monomial matrix $P$ over $\FF_q$ with 
$C' = C \cdot P = \{ x P \mid x \in C\}$.  
The {\em dual} code $C^{\perp}$ of an $[n,k]$ code $C$ over $\FF_q$
is defined as
$
C^{\perp}=
\{x \in \FF_q^n \mid x \cdot y = 0 \text{ for all } y \in C\},
$
where $x \cdot y$ is the standard inner product.

A code $C$ of length $n$ is called {\em linear complementary dual}
(or a linear code with complementary dual)
if $C \cap C^\perp = \{\0_n\}$, where $\0_n$ denotes the zero vector of length $n$.
We say that such a code is LCD for short.
LCD codes were introduced by Massey~\cite{Massey} and gave an optimum linear
coding solution for the two user binary adder channel.
Recently, much work has been done concerning LCD codes
for both theoretical and practical reasons (see~\cite{mf}, 
\cite{CMTQ2}, \cite{DKOSS}, \cite{bound}, \cite{HS}
and the references therein).
% It is a fundamental problem to classify LCD $[n,k]$ codes
% and determine the largest minimum weight
% among all LCD $[n,k]$ codes.
In particular, we emphasize the recent work by 
Carlet, Mesnager, Tang, Qi and Pellikaan~\cite{CMTQ2}.
It has been shown in~\cite{CMTQ2} that
any code over $\FF_q$ is equivalent to some LCD code
for $q \ge 4$.
This motivates us to study binary LCD codes.

From now on, all codes mean binary codes, and 
binary codes are simply called codes.
It is a fundamental problem to determine the largest minimum weight
among all LCD $[n,k]$ codes.
In this paper, we study the minimum weights of linear complementary 
dual codes.
Throughout this paper,
let $d(n,k)$ denote the largest minimum weight among all 
LCD $[n,k]$ codes.

It is trivial that $d(n,n)=1$.
It is known~\cite{DKOSS} that 
$d(n,1)=n$ and $n-1$  if $n$ is odd and even, respectively.
It was shown in~\cite{bound} that
$d(n,2) =
\lfloor \frac{2n}{3}\rfloor$ if  $n \equiv 1,2,3,4 \pmod{6}$, and
$\lfloor \frac{2n}{3} \rfloor-1$ otherwise for $n \ge 2$.
In addition, it was shown in~\cite{HS} that
$d(n,3) =\left\lfloor \frac{4n}{7}\right\rfloor$
if $n \equiv 3,5 \pmod {7}$ and 
$\left\lfloor \frac{4n}{7}\right\rfloor-1$ otherwise
for $n \ge 3$.
The aim of this paper is to establish the following theorems.

\begin{thm}\label{thm:dim4}
If $n \equiv 5,9,13 \pmod{15}$ and $n \ge 4$, then
\[
d(n,4) =
\left\lfloor \frac{8n}{15}\right\rfloor.
\]
If $n \equiv 2,3,4,6,10 \pmod{15}$ and $n \ge 4$, then
\[
d(n,4)=
\left\lfloor \frac{8n}{15}\right\rfloor-1.
\]
\end{thm}

\begin{thm}\label{thm:dim5}
If $n \equiv 3,5,7,11,19,20,22,26 \pmod{31}$ and $n \ge 5$, then
\[
d(n,5)=\left\lfloor \frac{16n}{31}\right\rfloor-1.
\]
If $n \equiv 4 \pmod{31}$ and $n \ge 5$, then
\[
d(n,5)=\left\lfloor \frac{16n}{31}\right\rfloor-2.
\]
\end{thm}

The values $d(n,k)$ were determined in~\cite{bound} and~\cite{HS}
for $n \le 12$ and $13 \le n \le 16$, respectively.
In this paper, we extend the results to lengths up to $24$.
To do this, 
we complete classifications of (unrestricted)  codes for
the parameters listed in Tables~\ref{Tab:New} and~\ref{Tab:New-2}.

All computer calculations in this paper were
done by programs in {\sc Magma}~\cite{Magma} and 
programs in the language C.
Two software libraries {\sc NTL}~\cite{NTL} and 
{\sc nauty and Traces}~\cite{nauty} were used.

%%%%%%%%%%%%%%%%%%%%%%%%%%%%%%%%%%%%%%%%%%%%%%
\section{LCD codes of dimension 4}
\label{Sec:k4}
Throughout this paper, we use the following notations.
Let $\0_{s}$ and $\1_{s}$ denote the zero vector and 
the all-one vector of length $s$, respectively.
Let $I_k$ denote the identity matrix of order $k$ and
let $A^T$ denote the transpose of a matrix $A$.
Let $\ZZ_{\ge 0}$ denote the set of nonnegative integers.

Suppose that there is an (unrestricted) $[n,4,d]$ code with $n \ge 4$.
By the Griesmer bound~\cite{Gr}, we have
\begin{equation}\label{eq:b4}
d \le
\begin{cases}
\left\lfloor \frac{8n}{15}\right\rfloor &\text{ if } 
n \equiv 0,1,5,7,8,9,11,12,13,14 \pmod{15},
\\
\left\lfloor \frac{8n}{15}\right\rfloor-1 &\text{ otherwise.}
\end{cases}
\end{equation}

For $a=(a_1,a_2,\ldots,a_{15})\in \ZZ_{\ge 0}^{15}$,
we define an $[n,4]$ code $C(a)$ having generator matrix of the form
$
G(a)=
\left(
\begin{array}{ccccc}
I_4 &  M(a)  \\
\end{array}
\right),
$
where 
\begin{multline*}
M(a)
=
\left(
\begin{array}{ccccccccccccccc}
\1_{a_{1}}&\1_{a_{2}}&\1_{a_{3}}&\1_{a_{4}}&\0_{a_{5}}&\1_{a_{6}}&\1_{a_{7}}\\
\1_{a_{1}}&\1_{a_{2}}&\1_{a_{3}}&\0_{a_{4}}&\1_{a_{5}}&\1_{a_{6}}&\0_{a_{7}}\\
\1_{a_{1}}&\1_{a_{2}}&\0_{a_{3}}&\1_{a_{4}}&\1_{a_{5}}&\0_{a_{6}}&\1_{a_{7}}\\
\1_{a_{1}}&\0_{a_{2}}&\1_{a_{3}}&\1_{a_{4}}&\1_{a_{5}}&\0_{a_{6}}&\0_{a_{7}}
\end{array}
\right.
\\
\left.
\begin{array}{ccccccccccccccc}
\1_{a_{8}}&\0_{a_{9}}&\0_{a_{10}}&\0_{a_{11}}&\1_{a_{12}}&\0_{a_{13}}&\0_{a_{14}}&\0_{a_{15}}\\
\0_{a_{8}}&\1_{a_{9}}&\1_{a_{10}}&\0_{a_{11}}&\0_{a_{12}}&\1_{a_{13}}&\0_{a_{14}}&\0_{a_{15}}\\
\0_{a_{8}}&\1_{a_{9}}&\0_{a_{10}}&\1_{a_{11}}&\0_{a_{12}}&\0_{a_{13}}&\1_{a_{14}}&\0_{a_{15}}\\
\1_{a_{8}}&\0_{a_{9}}&\1_{a_{10}}&\1_{a_{11}}&\0_{a_{12}}&\0_{a_{13}}&\0_{a_{14}}&\1_{a_{15}}
\end{array}
\right).
\end{multline*}
% We say that the minimum weight of $C^\perp$ 
% is the {\em dual distance} of $C$.
% Obviously, 
% the code $C(a)$ has dual distance at least $2$.
By considering all codewords, 
the weight enumerator of the code $C(a)$ is written using
$a_1,a_2,\ldots,a_{15}$ as follows:
\begin{align}\label{eq:WE}
&
1
+y^{1+a_1+a_2+a_3+a_4+a_6+a_7+a_8+a_{12}}
+y^{1+a_1+a_2+a_3+a_5+a_6+a_9+a_{10}+a_{13}}
\notag \\&
+y^{1+a_1+a_2+a_4+a_5+a_7+a_9+a_{11}+a_{14}}
+y^{1+a_1+a_3+a_4+a_5+a_8+a_{10}+a_{11}+a_{15}}
\notag \\&
+y^{2+a_4+a_5+a_7+a_8+a_9+a_{10}+a_{12}+a_{13}}
+y^{2+a_3+a_5+a_6+a_8+a_9+a_{11}+a_{12}+a_{14}}
\notag \\&
+y^{2+a_2+a_5+a_6+a_7+a_{10}+a_{11}+a_{12}+a_{15}}
+y^{2+a_3+a_4+a_6+a_7+a_{10}+a_{11}+a_{13}+a_{14}}
\\&
+y^{2+a_2+a_4+a_6+a_8+a_9+a_{11}+a_{13}+a_{15}}
+y^{2+a_2+a_3+a_7+a_8+a_9+a_{10}+a_{14}+a_{15}}
\notag \\&
+y^{3+a_1+a_2+a_8+a_{10}+a_{11}+a_{12}+a_{13}+a_{14}}
+y^{3+a_1+a_3+a_7+a_9+a_{11}+a_{12}+a_{13}+a_{15}}
\notag \\&
+y^{3+a_1+a_4+a_6+a_9+a_{10}+a_{12}+a_{14}+a_{15}}
+y^{3+a_1+a_5+a_6+a_7+a_8+a_{13}+a_{14}+a_{15}}
\notag \\&
+y^{4+a_2+a_3+a_4+a_5+a_{12}+a_{13}+a_{14}+a_{15}}.
\notag
\end{align}
The $(i,j)$-entry $b_{i,j}$ of $G(a)G(a)^T$ is written using
$a_1,a_2,\ldots,a_{15}$ as follows:
\begin{align}
b_{1,1}=&1+a_1+a_2+a_3+a_4+a_6+a_7+a_8+a_{12},\notag \\
b_{1,2}=&a_1+a_2+a_3+a_6,\notag \\
b_{1,3}=&a_1+a_2+a_4+a_7,\notag \\
b_{1,4}=&a_1+a_3+a_4+a_8,\notag \\
b_{2,2}=&1+a_1+a_2+a_3+a_5+a_6+a_9+a_{10}+a_{13},\label{eq:det}\\
b_{2,3}=&a_1+a_2+a_5+a_9,\notag \\
b_{2,4}=&a_1+a_3+a_5+a_{10},\notag \\
b_{3,3}=&1+a_1+a_2+a_4+a_5+a_7+a_9+a_{11}+a_{14},\notag \\
b_{3,4}=&a_1+a_4+a_5+a_{11},\notag \\
b_{4,4}=&1+a_1+a_3+a_4+a_5+a_8+a_{10}+a_{11}+a_{15}.\notag 
\end{align}
% Thus, the determinant $\det(G(a)G(a)^T)$ is written using
% $a_1,a_2,\ldots,a_{15}$.

%%%%%%%%%%%%%%%%%%%%%%%%%%%%%%
\begin{table}[thb]
\caption{LCD codes of dimension 4}
\label{Tab:dim4-1}
\begin{center}
%{\small
{\footnotesize
%{\scriptsize
%{\tiny
\begin{tabular}{c|l}
\noalign{\hrule height0.8pt}
Code & \multicolumn{1}{c}{$a$} \\
\hline
$C_{15t+2}$ & $(t,t,t,t,t,t,t,t-1,t+1,t+1,t-1,t,t-1,t-1,t)$\\
$C_{15t+3}$ & $(t,t,t,t,t,t,t,t,t,t,t-1,t,t,t,t)$\\
$C_{15t+4}$ & $(t,t,t,t,t,t,t,t,t,t,t,t,t,t,t)$\\
$C_{15t+5}$ & $(t,t,t,t,t,t,t,t+1,t+1,t,t+1,t,t,t-1,t-1)$\\
$C_{15t+6}$ & $(t,t,t,t,t,t,t,t+1,t+1,t,t,t,t,t,t)$\\
$C_{15t+9}$ & $(t,t,t,t,t+1,t+1,t+1,t+1,t,t+1,t+1,t,t,t,t-1)$\\
$C_{15t+10}$ &$(t,t,t,t,t,t,t+1,t+1,t+1,t+2,t+1,t+1,t,t,t-1)$\\
$C_{15t+13}$ &$(t,t,t,t+1,t+1,t,t+1,t+2,t+2,t+1,t+1,t+1,t+1,t-1,t-1)$\\
\hline
$C_{15t+1}$ & $(t,t,t,t,t,t,t,t,t,t-1,t-1,t,t,t-1,t)$\\
$C_{15t+7}$ & $(t,t,t,t,t,t,t,t+1,t+1,t,t,t,t,t,t+1)$\\
$C_{15t+8}$ & $(t,t,t,t,t,t,t+1,t+1,t+1,t+1,t,t,t,t,t)$\\
$C_{15t+11}$ &$(t,t,t,t,t,t,t,t+1,t+1,t+1,t+2,t+2,t+1,t,t-1)$\\
$C_{15t+12}$ &$(t,t,t,t,t+1,t,t+1,t+2,t+2,t+1,t+1,t+1,t+1,t-1,t-1)$\\
$C_{15t+14}$ &$(t,t,t,t,t+1,t,t+1,t+2,t+2,t+2,t+1,t+2,t,t,t-1)$\\
\hline
$C_{15t}$ & $(t,t,t,t,t,t,t,t,t,t-1,t-1,t,t,t-1,t-1)$\\
\noalign{\hrule height0.8pt}
\end{tabular}
}
\end{center}
\end{table}
%%%%%%%%%%%%%%%%%%%%%%%%%%%%%%%

Write $n=15t+s$, where $t$ is a nonnegative integer
and $s\in\{0,1,\ldots,14\}$.
For $s=2,3,4,5,6,9,10,13$, 
by considering $C(a)$, we found the codes
$C_{15t+s}$ 
meeting the bound~\eqref{eq:b4} with equality,
where the vectors $a$ are listed in Table~\ref{Tab:dim4-1}.
The minimum weights are determined from 
the weight enumerators $W$ obtained by~\eqref{eq:WE},
where $W$ are listed in Table~\ref{Tab:dim4-WE}.

It was shown in~\cite{Massey} that a code $C$ is LCD if and only if 
$G G^T$ is nonsingular for any generator matrix $G$ of $C$.
This fact is used in order to show that a given code is LCD,
throughout this paper.
From~\eqref{eq:det}, 
the determinants $\det(G(a)G(a)^T)$ are written using $t$,
where the results are listed in Table~\ref{Tab:dim4-det}.
It follows from the table that 
$\det(G(a)G(a)^T)=1$ for every nonnegative integer $t$.
Hence,  we have the following:
\begin{itemize}
\item[\rm (1)]$C_{15t+ 2}$ is an LCD $[15t+ 2,4,  8t]$ code $(t \ge 1)$,
\item[\rm (2)]$C_{15t+ 3}$ is an LCD $[15t+ 3,4,  8t]$ code $(t \ge 1)$,
\item[\rm (3)]$C_{15t+ 4}$ is an LCD $[15t+ 4,4,8t+1]$ code $(t \ge 0)$,
\item[\rm (4)]$C_{15t+ 5}$ is an LCD $[15t+ 5,4,8t+2]$ code $(t \ge 1)$,
\item[\rm (5)]$C_{15t+ 6}$ is an LCD $[15t+ 6,4,8t+2]$ code $(t \ge 0)$,
\item[\rm (6)]$C_{15t+ 9}$ is an LCD $[15t+ 9,4,8t+4]$ code $(t \ge 1)$,
\item[\rm (7)]$C_{15t+10}$ is an LCD $[15t+10,4,8t+4]$ code $(t \ge 1)$,
\item[\rm (8)]$C_{15t+13}$ is an LCD $[15t+13,4,8t+6]$ code $(t \ge 1)$.
\end{itemize}
% These codes have dual distance $2$.
In addition, there is an LCD $[n,4,d]$ code for
$(n,d)=(5,2),( 9,4),(10,4)$ and $(13,6)$~\cite{bound} and~\cite{HS}.
This completes the proof of Theorem~\ref{thm:dim4}.

%%%%%%%%%%%%%%%%%%%%%%%%%%%%%%
\begin{table}[thb]
\caption{Weight enumerators of $C_{15t+s}$}
\label{Tab:dim4-WE}
\begin{center}
%{\small
{\footnotesize
%{\scriptsize
%{\tiny
\begin{tabular}{c|l}
\noalign{\hrule height0.8pt}
Code & \multicolumn{1}{c}{$W$} \\
\hline
$C_{15t+2}$ &$1+8y^{8t}+6y^{8t+2}+y^{8t+4}$ \\
$C_{15t+3}$ &$1+2y^{8t}+6y^{8t+1}+4y^{8t+2}+2y^{8t+3}+y^{8t+4}$ \\
$C_{15t+4}$ &$1+4y^{8t+1}+6y^{8t+2}+4y^{8t+3}+y^{8t+4}$ \\
$C_{15t+5}$ &$1+10y^{8t+2}+5y^{8t+4}$ \\
$C_{15t+6}$ &$1+6y^{8t+2}+9y^{8t+4}$ \\
$C_{15t+9}$ &$1+9y^{8t+4}+6y^{8t+6}$ \\
$C_{15t+10}$ &$1+7y^{8t+4}+6y^{8t+6}+2y^{8t+8}$ \\
$C_{15t+13}$ &$1+10y^{8t+6}+4y^{8t+8}+y^{8t+12}$ \\
\hline
$C_{15t+1}$ &$1+3y^{8t-1}+5y^{8t}+4y^{8t+1}+2y^{8t+2}+y^{8t+3}$ \\
$C_{15t+7}$ &$1+4y^{8t+2}+2y^{8t+3}+3y^{8t+4}+6y^{8t+5}$ \\
$C_{15t+8}$ &$1+4y^{8t+3}+5y^{8t+4}+4y^{8t+5}+2y^{8t+6}$ \\
$C_{15t+11}$ &$1+6y^{8t+4}+5y^{8t+6}+3y^{8t+8}+y^{8t+10}$ \\
$C_{15t+12}$ &$1+6y^{8t+5}+4y^{8t+6}+y^{8t+7}+3y^{8t+8}+y^{8t+11}$ \\
$C_{15t+14}$ &$1+8y^{8t+6}+4y^{8t+8}+2y^{8t+10}+y^{8t+12}$ \\
\hline
$C_{15t}$ &$1+y^{8t-2}+4y^{8t-1}+5y^{8t}+4y^{8t+1}+y^{8t+2}$ \\
\noalign{\hrule height0.8pt}
\end{tabular}
}
\end{center}
\end{table}
%%%%%%%%%%%%%%%%%%%%%%%%%%%%%%%

%%%%%%%%%%%%%%%%%%%%%%%%%%%%%%
\begin{table}[thb]
\caption{$\det(G(a)G(a)^T)$ for $C_{15t+s}$}
\label{Tab:dim4-det}
\begin{center}
%{\small
{\footnotesize
%{\scriptsize
%{\tiny
\begin{tabular}{c|l}
\noalign{\hrule height0.8pt}
Code & \multicolumn{1}{c}{$\det(G(a)G(a)^T)$} \\
\hline
$C_{15t+2}$ &$1280t^4+512t^3-96t^2-32t+1$\\
$C_{15t+3}$ &$1280t^4+640t^3+64t^2-8t-1$\\
$C_{15t+4}$ &$1280t^4+1024t^3+288t^2+32t+1$\\
$C_{15t+5}$ &$1280t^4+1664t^3+688t^2+104t+5$\\
$C_{15t+6}$ &$1280t^4+1792t^3+800t^2+144t+9$\\
$C_{15t+9}$ &$1280t^4+3072t^3+2656t^2+976t+129$\\
$C_{15t+10}$ &$1280t^4+3328t^3+3040t^2+1152t+153$\\
$C_{15t+13}$ &$1280t^4+4480t^3+5424t^2+2744t+493$\\
\hline
$C_{15t+1}$ &$1280t^4-80t^2+1$\\
$C_{15t+7}$ &$1280t^4+2048t^3+1056t^2+216t+15$\\
$C_{15t+8}$ &$1280t^4+2560t^3+1760t^2+480t+45$\\
$C_{15t+11}$ &$1280t^4+3456t^3+3248t^2+1240t+161$\\
$C_{15t+12}$ &$1280t^4+4096t^3+4496t^2+2064t+339$\\
$C_{15t+14}$ &$1280t^4+4736t^3+6064t^2+3208t+597$\\
\hline
$C_{15t}$ &$1280t^4-256t^3-80t^2+8t+1$\\
\noalign{\hrule height0.8pt}
\end{tabular}
}
\end{center}
\end{table}
%%%%%%%%%%%%%%%%%%%%%%%%%%%%%%%

For $s=0,1,7,8,11,12$, 
similarly, we found the codes $C_{15t+s}=C(a)$,
where the vectors $a$ are listed
in Table~\ref{Tab:dim4-1}.
These codes have  minimum weight one or two less than
the largest possible minimum weight in the bound~\eqref{eq:b4}.
Their weight enumerators $W$ obtained by~\eqref{eq:WE} and
their determinants $\det(G(a)G(a)^T)$ obtained by~\eqref{eq:det}
are listed
in Tables~\ref{Tab:dim4-WE} and~\ref{Tab:dim4-det},
respectively.
Hence, we have the following:
\begin{itemize}
\item[\rm (1)]$C_{15t   }$ is an LCD $[15t,   4,8t-2]$ code $(t \ge 1)$,
\item[\rm (2)]$C_{15t+ 1}$ is an LCD $[15t+ 1,4,8t-1]$ code $(t \ge 1)$,
\item[\rm (3)]$C_{15t+ 7}$ is an LCD $[15t+ 7,4,8t+2]$ code $(t \ge 0)$,
\item[\rm (4)]$C_{15t+ 8}$ is an LCD $[15t+ 8,4,8t+3]$ code $(t \ge 0)$,
\item[\rm (5)]$C_{15t+11}$ is an LCD $[15t+11,4,8t+4]$ code $(t \ge 1)$,
\item[\rm (6)]$C_{15t+12}$ is an LCD $[15t+12,4,8t+5]$ code $(t \ge 1)$,
\item[\rm (7)]$C_{15t+14}$ is an LCD $[15t+14,4,8t+6]$ code $(t \ge 1)$.
\end{itemize}
In addition, there is an LCD $[n,4,d]$ code for
$(n,d)=(11,4),(12,5)$ and $(14,6)$~\cite{bound} and~\cite{HS}.
Therefore, 
% there is an LCD $[n,4,\lfloor \frac{8n}{15}\rfloor-1]$ code
% for $n \equiv 1,7,8,11,12,14 \pmod{15}$ and $n \ge 4$, and
% there is an LCD $[n,4,\lfloor \frac{8n}{15}\rfloor-2]$ code
% for $n \equiv 0 \pmod{15}$ and $n \ge 4$.
we have the following:

\begin{prop}
If $n \equiv 1,7,8,11,12,14 \pmod{15}$ and $n \ge 4$, then
\[
d(n,4) =
\left\lfloor \frac{8n}{15}\right\rfloor \text{ or }
\left\lfloor \frac{8n}{15}\right\rfloor -1.
\]
If $n \equiv 0 \pmod{15}$ and $n \ge 4$, then
\[
d(n,4) =
\left\lfloor \frac{8n}{15}\right\rfloor,
\left\lfloor \frac{8n}{15}\right\rfloor -1  \text{ or }
\left\lfloor \frac{8n}{15}\right\rfloor -2.
\]
\end{prop}

Now we are a position to consider the existence of 
an LCD $[n,4,\lfloor \frac{8n}{15}\rfloor]$ code
for $n \equiv 1,7,8,11,12 \pmod{15}$  and
an LCD $[n,4,\lfloor \frac{8n}{15}\rfloor-1]$ code
for $n \equiv 0 \pmod{15}$.
There is no such code for $n \le 16$~\cite{bound} and~\cite{HS}.
For 
\begin{equation}\label{eq:dim4}
(n,d)=(22, 11),
(23, 12),
(26, 13),
(27, 14),
(30, 16),
(30, 15),
(31, 16),
\end{equation}
in order to verify that there is no LCD $[n,4,d]$
code,
our computer calculation completed a classification of 
(unrestricted) $[n,4,d]$ codes
by using the following method.
A {\em shortened code} $C'$ of a code $C$ is the set of all codewords 
in $C$ which are $0$ in a fixed coordinate with that 
coordinate deleted.
A shortened code $C'$ of an $[n,k,d]$ code $C$ with $d \ge 2$
is an $[n-1,k,d]$ code if the deleted coordinate
is a zero coordinate and an $[n-1,k-1,d']$ 
code with $d' \ge d$ 
% and covering radius $R\ge d-1$ 
otherwise.
An $[n,k,d]$ code $C$ gives $n$ shortened codes
and at least $k$ codes among them are $[n-1,k-1,d']$ codes
with $d' \ge d$.
Hence, 
by considering the inverse operation of shortening, 
any $[n,k,d]$ code with $d \ge 2$ is constructed from some
$[n-1,k-1,d']$ code with $d' \ge d$.
% and covering radius $R\ge d-1$.
% $R\ge d-1$ as follows.
% Let $C'$ be an $[n-1,k-1,d']$ code with $d' \ge d$.
% Up to equivalence, 
% we may assume that $C'$ has 
% a generator matrix of the form $
% \left(\begin{array}{cc}
% I_{k-1} & A 
% \end{array}\right)$.
% Then, up to equivalence, 
% an $[n,k,d]$ code, which is constructed from $C'$
% by considering the inverse operation of shortening, 
% has the following generator matrix
% \begin{equation}\label{eq:S1}
% \left(\begin{array}{c|ccc|ccccccc}
% 1&0& \cdots& 0 &b_1 &\cdots&b_{n-k}\\
% \hline
% 0& &       &&      & &  &\\
% \vdots& &I_{k-1}& & &A &\\
% 0& &       &      & &  &\\
% \end{array}\right),
% \end{equation}
% where $b=(b_1,b_2,\ldots,b_{n-k}) \in\FF_2^{n-k}$
% with $\wt(b) \ge d-1$.
In order to illustrate this method, we describe how 
$[22, 4, 11]$ codes were classified.
Let $d_{\text{all}}(n,k)$ denote the largest minimum weight
among all (unrestricted) $[n,k]$ codes.
From~\cite{G}, we know
\[
d_{\text{all}}(21,3)=12 \text{ and }
d_{\text{all}}(20,2)=13.
\]
We first classified all (unrestricted) $[20,2,d]$ codes with
$d=11,12,13$ by a direct method.
From this classification,
by the above method, we found all inequivalent
$[21,3,d]$ codes with $d=11,12$.
Then we found  all inequivalent $[22,4,11]$ codes.

Let $N_{n,k,d}$ denote the number of all inequivalent $[n,k,d]$ codes.
In Table~\ref{Tab:dim4-small},
we list $N_{n,4,d}$ and 
$N_{n,k,d'}$ ($k=2,3, d'=d,d+1,d+2,d+3$) 
for $(n,d)$ in~\eqref{eq:dim4}.
In order to give generator matrices 
$
\left(
\begin{array}{ccccc}
I_4 &  M_{n,i}  \\
\end{array}
\right)
$
of all inequivalent $[n,4,d]$ codes, 
we only list the four rows $m_1,m_2,m_3,m_4$
of $M_{n,i}$ in Table~\ref{Tab:dim4-GM}.
To save space, the sequences $m=(m_1,m_2,m_3,m_4)$
are written in octal using $0=(000),$ $1=(001),\ldots,\ 7=(111)$,
together with $a=(0)$ and $b=(1)$.  
Note that 
$M_{30,1}$ corresponds to $d=16$ and
$M_{30,i}$ $(i=2,\ldots,10)$ correspond to $d=15$.

%%%%%%%%%%%%%%%%%%%%%%%%%%%%%%
\begin{table}[thbp]
\caption{Numbers of $[n,4,d]$ codes for $(n,d)$ in~\eqref{eq:dim4}}
\label{Tab:dim4-small} 
\begin{center}
%{\small
{\footnotesize
%{\scriptsize
%{\tiny
\begin{tabular}{l|ccccccc}
\noalign{\hrule height0.8pt}
\multicolumn{1}{c|}{$(n,d)$ } 
& $(22, 11)$ & $(23, 12)$ & $(26, 13)$&$(27,14)$&$(30,16)$&$(30,15)$
&$(31,16)$\\
\hline
$N_{n,4,d}$     &2 &1 &2 &1 &1 &9 &5\\
\hline
$N_{n-1,3,d}$   &6 &4 &13&7 &4 &27&16\\
$N_{n-1,3,d+1}$ &1 &- &1 &- &- &4 &-\\
\hline
$N_{n-2,2,d}$   &10&10&15&15&15&21&21\\
$N_{n-2,2,d+1}$ &6 &3 &10&6 &6 &15&10\\
$N_{n-2,2,d+2}$ &1 &1 &3 &3 &3 &6 &6 \\
$N_{n-2,2,d+3}$ &- &- &1 &- &- &3 &1 \\
\noalign{\hrule height0.8pt}
\end{tabular}
}
\end{center}
\end{table}
%%%%%%%%%%%%%%%%%%%%%%%%%%%%%%%

%%%%%%%%%%%%%%%%%%%%%%%%%%%%%%
\begin{table}[thb]
\caption{Generator matrices of $[n,4,d]$ codes for $(n,d)$ in~\eqref{eq:dim4}}
\label{Tab:dim4-GM} 
\begin{center}
%{\small
{\footnotesize
%{\scriptsize
%{\tiny
\begin{tabular}{c|l}
\noalign{\hrule height0.8pt}
$M_{22,1}$&617170773600001777475345\\
$M_{22,2}$&633330767460001777475345\\
$M_{23,1}$&7066743767400003777533415b\\
$M_{26,1}$&74607433743630000077774773714a\\
$M_{26,2}$&63653061761714000077771676540b\\
$M_{27,1}$&760663616177700000177775750755ba\\
$M_{30,1}$&7074633617703754000007777766174433ab\\
$M_{30,2}$&3746066317361730000003777767251176ab\\
$M_{30,3}$&5147543306363674000003777767251176ab\\
$M_{30,4}$&7436630317741714000003777751676754ba\\
$M_{30,5}$&7306663617773600000003777764564745aa\\
$M_{30,6}$&3615263617767700000003777764564745aa\\
$M_{30,7}$&7314633617743740000003777764564745aa\\
$M_{30,8}$&7707043617363614000003777764564745aa\\
$M_{30,9}$&7317063617761714000003777764564745aa\\
$M_{30,10}$&3615700757303746000003777764564745aa\\
$M_{31,1}$&554633154717077600000077777672511762\\
$M_{31,2}$&730661730777077000000077777137753663\\
$M_{31,3}$&347474036774170660000077777516767544\\
$M_{31,4}$&760374630777633000000077777172511762\\
$M_{31,5}$&547433154774077600000077777172511762\\
\noalign{\hrule height0.8pt}
\end{tabular}
}
\end{center}
\end{table}

From the above classification,
our computer calculation shows the following result.
 
\begin{prop}
There is no LCD $[n,4,d]$ code for 
\[
(n,d)=(22, 11),
(23, 12),
(26, 13),
(27, 14),
(30, 16),
(30, 15),
(31, 16).
\]
\end{prop}

%%%%%%%%%%%%%%%%%%%%%%%%%%%%%%%%%%%%%%%%%%%%%%
\section{LCD codes of dimension 5}
\label{Sec:k5}
Suppose that there is an (unrestricted) $[n,5,d]$ code with $n \ge 5$.
By the Griesmer bound~\cite{Gr}, we have
\begin{equation}\label{eq:b5}
d \le
\begin{cases}
\left\lfloor \frac{16n}{31}\right\rfloor &\text{ if } 
n \equiv 0,1,9,13,15,16, \\
& \hspace*{1.5cm} 17,21,23,24,25,27,28,29,30 \pmod{31}, \\
\left\lfloor \frac{16n}{31}\right\rfloor-1 &\text{ if } 
n \equiv 2,3,5,6,7,8,10,\\
& \hspace*{1.5cm} 11,14,18,19,20,22,26 \pmod{31}, \\
\left\lfloor \frac{16n}{31}\right\rfloor-2 &\text{ otherwise. } 
\end{cases}
\end{equation}

For $a=(a_1,a_2,\ldots,a_{31})\in \ZZ_{\ge 0}^{31}$,
we define an $[n,5]$ code $C(a)$ having generator matrix 
of the form
$
G(a)=
\left(
\begin{array}{cccccc}
I_5 &  M(a)  \\
\end{array}
\right),
$
where $M(a)$ is listed in Figure~\ref{Fig}.
Using an approach similar to that in the previous section,
the weight enumerator $W$ and the determinant $\det(G(a)G(a)^T)$
for the code $C(a)$ are written using
$a_1,a_2,\ldots,a_{31}$.
% To save space, we do not list the details.

%%%%%%%%%%%%%%%%%%%%%%%%%%%%%%%%%%%%%%%%%%%%%%%%%%%%
\begin{figure}[thb]
\label{Fig}
\begin{center}
%{\small
{\footnotesize
%{\scriptsize
%{\tiny
\begin{align*}
&
\left(
\begin{array}{ccccccccccccccc}
\1_{a_{1}}&\1_{a_{2}}&\1_{a_{3}}&\1_{a_{4}}&\1_{a_{5}}&\1_{a_{6}}&\1_{a_{7}}\\
\1_{a_{1}}&\1_{a_{2}}&\1_{a_{3}}&\1_{a_{4}}&\1_{a_{5}}&\1_{a_{6}}&\1_{a_{7}}\\
\1_{a_{1}}&\1_{a_{2}}&\1_{a_{3}}&\1_{a_{4}}&\0_{a_{5}}&\0_{a_{6}}&\0_{a_{7}}\\
\1_{a_{1}}&\1_{a_{2}}&\0_{a_{3}}&\0_{a_{4}}&\1_{a_{5}}&\1_{a_{6}}&\0_{a_{7}}\\
\1_{a_{1}}&\0_{a_{2}}&\1_{a_{3}}&\0_{a_{4}}&\1_{a_{5}}&\0_{a_{6}}&\1_{a_{7}}
\end{array}
\right.
\\
&\qquad
\left.
\begin{array}{ccccccccccccccc}
\1_{a_{8}} &\1_{a_{9}} &\1_{a_{10}}&\1_{a_{11}}&\1_{a_{12}}&\1_{a_{13}}&\1_{a_{14}}&\1_{a_{15}}\\
\1_{a_{8}} &\0_{a_{9}} &\0_{a_{10}}&\0_{a_{11}}&\0_{a_{12}}&\0_{a_{13}}&\0_{a_{14}}&\0_{a_{15}}\\
\0_{a_{8}} &\1_{a_{9}} &\1_{a_{10}}&\1_{a_{11}}&\1_{a_{12}}&\0_{a_{13}}&\0_{a_{14}}&\0_{a_{15}}\\
\0_{a_{8}} &\1_{a_{9}} &\1_{a_{10}}&\0_{a_{11}}&\0_{a_{12}}&\1_{a_{13}}&\1_{a_{14}}&\0_{a_{15}}\\
\0_{a_{8}} &\1_{a_{9}} &\0_{a_{10}}&\1_{a_{11}}&\0_{a_{12}}&\1_{a_{13}}&\0_{a_{14}}&\1_{a_{15}}
\end{array}
\right.
\\
&\qquad
\left.
\begin{array}{ccccccccccccccc}
\1_{a_{16}}&\0_{a_{17}}&\0_{a_{18}}&\0_{a_{19}}&\0_{a_{20}}&\0_{a_{21}}&\0_{a_{22}}&\0_{a_{23}}\\
\0_{a_{16}}&\1_{a_{17}}&\1_{a_{18}}&\1_{a_{19}}&\1_{a_{20}}&\1_{a_{21}}&\1_{a_{22}}&\1_{a_{23}}\\
\0_{a_{16}}&\1_{a_{17}}&\1_{a_{18}}&\1_{a_{19}}&\1_{a_{20}}&\0_{a_{21}}&\0_{a_{22}}&\0_{a_{23}}\\
\0_{a_{16}}&\1_{a_{17}}&\1_{a_{18}}&\0_{a_{19}}&\0_{a_{20}}&\1_{a_{21}}&\1_{a_{22}}&\0_{a_{23}}\\
\0_{a_{16}}&\1_{a_{17}}&\0_{a_{18}}&\1_{a_{19}}&\0_{a_{20}}&\1_{a_{21}}&\0_{a_{22}}&\1_{a_{23}}
\end{array}
\right.
\\
&\qquad
\left.
\begin{array}{ccccccccccccccc}
\0_{a_{24}}&\0_{a_{25}}&\0_{a_{26}}&\0_{a_{27}}&\0_{a_{28}}&\0_{a_{29}}&\0_{a_{30}}&\0_{a_{31}}\\
\1_{a_{24}}&\0_{a_{25}}&\0_{a_{26}}&\0_{a_{27}}&\0_{a_{28}}&\0_{a_{29}}&\0_{a_{30}}&\0_{a_{31}}\\
\0_{a_{24}}&\1_{a_{25}}&\1_{a_{26}}&\1_{a_{27}}&\1_{a_{28}}&\0_{a_{29}}&\0_{a_{30}}&\0_{a_{31}}\\
\0_{a_{24}}&\1_{a_{25}}&\1_{a_{26}}&\0_{a_{27}}&\0_{a_{28}}&\1_{a_{29}}&\1_{a_{30}}&\0_{a_{31}}\\
\0_{a_{24}}&\1_{a_{25}}&\0_{a_{26}}&\1_{a_{27}}&\0_{a_{28}}&\1_{a_{29}}&\0_{a_{30}}&\1_{a_{31}}
\end{array}
\right)
\end{align*}
}
\end{center}
\caption{Matrix $M(a)$}
\end{figure}
%%%%%%%%%%%%%%%%%%%%%%%%%%%%%%%%%%%%%%%%%%%%%%%%%%%%

Write $n=31t+s$, where $t$ is a nonnegative integer
and $s\in\{0,1,\ldots,30\}$.
For $s=
 3,
 4,
 5,
 7,
11,
19,
20,
22,
26$,
by considering $C(a)$, we found the codes
$D_{31t+s}$ 
meeting the bound~\eqref{eq:b5} with equality,
where the vectors $a$ are listed in Table~\ref{Tab:dim5}.
In Table~\ref{Tab:dim5}, 
we denote $t-1,t+1$ by $t_-,t_+$, respectively.
The minimum weights are determined from 
the weight enumerators $W$,
where $W$ are listed in Table~\ref{Tab:dim5-WE}.

The determinants $\det(G(a)G(a)^T)$ are written using $t$,
where the results are listed in Table~\ref{Tab:dim5-det}.
It follows from the table that 
$\det(G(a)G(a)^T)=1$ for every nonnegative integer $t$.
Hence, we have the following:
\begin{itemize}
\item[\rm (1)]
$D_{31t+ 3}$ is an LCD $[31t + 3, 5,16t]$ code $(t \ge 1)$,
\item[\rm (2)]
$D_{31t+ 4}$ is an LCD $[31t + 4, 5,16t]$ code $(t \ge 1)$,
\item[\rm (3)]
$D_{31t+ 5}$ is an LCD $[31t + 5, 5,16t + 1]$ code  $(t \ge 0)$,
\item[\rm (4)]
$D_{31t+ 7}$ is an LCD $[31t + 7, 5,16t + 2]$ code $(t \ge 1)$,
\item[\rm (5)]
$D_{31t+11}$ is an LCD $[31t +11, 5,16t + 4]$ code $(t \ge 1)$,
\item[\rm (6)]
$D_{31t+19}$ is an LCD $[31t +19, 5,16t + 8]$ code $(t \ge 1)$,
\item[\rm (7)]
$D_{31t+20}$ is an LCD $[31t +20, 5,16t + 9]$ code $(t \ge 1)$,
\item[\rm (8)]
$D_{31t+22}$ is an LCD $[31t +22, 5,16t + 10]$ code $(t \ge 1)$,
\item[\rm (9)]
$D_{31t+26}$ is an LCD $[31t +26, 5,16t + 12]$ code $(t \ge 1)$.
\end{itemize}
% These codes have dual distance $2$.
In addition, 
there is an LCD $[n,5,d]$ code
for $(n,d)=(7,2)$ and $(11,4)$~\cite{bound} and~\cite{HS}.
Our computer search found an LCD $[n,5,d]$ code
%  with dual distance $2$
for $(n,5)=(19,8),(20,9),(22,10)$ and $(26,12)$.
These codes have generator matrices
$
\left(
\begin{array}{cc}
 I_5  & M_i \\
\end{array}
\right),
$
% where $M_{i}$ $(i=19,20,22,26)$ are listed in Figure~\ref{Fig:2},
% respectively.
where the five rows of 
$M_{i}$ $(i=19,20,22,26)$ are listed in 
Table~\ref{Tab:dim4-GM0}, respectively.
% These codes have the following weight enumerators:
% \begin{align*}
% &
% 1+7y^8+10y^9+5y^{10}+5y^{11}+2y^{12}+y^{14}+y^{15},
% \\&
% 1+ 10y^9 + 10y^{10}+ 5y^{11}+ 5y^{12}+ y^{15},
% \\&
% 1+10y^{10}+10y^{11}+ 5y^{12}+ 5y^{13}+ y^{17},
% \\&
% 1+ 9y^{12}+ 9y^{13}+ 6y^{14}+ 6y^{15}+ y^{17},
% \end{align*}
% respectively.
This completes the proof of Theorem~\ref{thm:dim5}.

%%%%%%%%%%%%%%%%%%%%%%%%%%%%%%
\begin{table}[thbp]
\caption{Weight enumerators of $D_{31t+s}$}
\label{Tab:dim5-WE}
\begin{center}
%{\small
{\footnotesize
%{\scriptsize
%{\tiny
\begin{tabular}{c|l}
\noalign{\hrule height0.8pt}
Code & \multicolumn{1}{c}{$W$} \\
\hline
$D_{31t+3}$ & 
$1+8y^{16t}+9y^{16t+1}+6y^{16t+2}+6y^{16t+3}+y^{16t+4}+y^{16t+5}$\\
$D_{31t+4}$ &
$1+2y^{16t}+9y^{16t+1}+10y^{16t+2}+6y^{16t+3}+3y^{16t+4}+y^{16t+5}$\\
$D_{31t+5}$ &
$1+5y^{16t+1}+10y^{16t+2}+10y^{16t+3}+5y^{16t+4}+y^{16t+5}$\\
$D_{31t+7}$ &
$1+6y^{16t+2}+9y^{16t+3}+9y^{16t+4}+6y^{16t+5}+y^{16t+7}$\\
$D_{31t+11}$&
$1+7y^{16t+4}+9y^{16t+5}+6y^{16t+6}+6y^{16t+7}+2y^{16t+8}+y^{16t+9}$\\
$D_{31t+19}$&
$1+8y^{16t+8}+8y^{16t+9}+6y^{16t+10}+6y^{16t+11}+y^{16t+12}+y^{16t+13}
+y^{16t+17}$ \\
$D_{31t+20}$&
$1+10y^{16t+9}+10y^{16t+10}+5y^{16t+11}+5y^{16t+12}+y^{16t+15}$\\
$D_{31t+22}$&
$1+10y^{16t+10}+10y^{16t+11}+5y^{16t+12}+5y^{16t+13}+y^{16t+17}$\\
$D_{31t+26}$&
$1+9y^{16t+12}+9y^{16t+13}+6y^{16t+14}+6y^{16t+15}+y^{16t+17}$\\
\hline
$D_{31t+1}$ &
$1+9y^{16t-1}+8y^{16t}+6y^{16t+1}+6y^{16t+2}+y^{16t+3}+y^{16t+4}$\\
$D_{31t+2}$ &
$1+3y^{16t-1}+8y^{16t}+10y^{16t+1}+6y^{16t+2}+3y^{16t+3}+y^{16t+4}$\\
$D_{31t+6}$ &
$1+3y^{16t+1}+6y^{16t+2}+10y^{16t+3}+9y^{16t+4}+3y^{16t+5}$\\
$D_{31t+8}$ &
$1+4y^{16t+2}+9y^{16t+3}+7y^{16t+4}+6y^{16t+5}+2y^{16t+6}+y^{16t+7}+2y^{16t+8}$\\
$D_{31t+9}$ &
$1+6y^{16t+3}+9y^{16t+4}+9y^{16t+5}+6y^{16t+6}+y^{16t+9}$\\
$D_{31t+10}$ &
$1+6y^{16t+3}+8y^{16t+4}+5y^{16t+5}+5y^{16t+6}+4y^{16t+7}+y^{16t+8}$\\
& $+y^{16t+9}+y^{16t+10}$\\
$D_{31t+12}$ &
$1+6y^{16t+4}+8y^{16t+5}+5y^{16t+6}+5y^{16t+7}+3y^{16t+8}+2y^{16t+9}$\\
&$+y^{16t+10}+y^{16t+11}$\\
$D_{31t+13}$ &
$1+8y^{16t+5}+9y^{16t+6}+5y^{16t+7}+5y^{16t+8}+2y^{16t+9}+y^{16t+10}+y^{16t+11}$\\
$D_{31t+14}$ &
$1+4y^{16t+5}+5y^{16t+6}+9y^{16t+7}+9y^{16t+8}+2y^{16t+9}+y^{16t+10}+y^{16t+11}$\\
$D_{31t+15}$ &
$1+8y^{16t+6}+9y^{16t+7}+4y^{16t+8}+6y^{16t+9}+2y^{16t+10}+y^{16t+11}+y^{16t+12}$\\
$D_{31t+17}$ &
$1+9y^{16t+7}+8y^{16t+8}+4y^{16t+9}+6y^{16t+10}+y^{16t+11}+y^{16t+12}+2y^{16t+13}$\\
$D_{31t+18}$ &
$1+6y^{16t+7}+7y^{16t+8}+6y^{16t+9}+5y^{16t+10}+3y^{16t+11}+2y^{16t+12}$\\
&$+y^{16t+14}+y^{16t+15}$\\
$D_{31t+21}$ &
$1+6y^{16t+9}+6y^{16t+10}+9y^{16t+11}+9y^{16t+12}+y^{16t+15}$\\
$D_{31t+23}$ &
$1+4y^{16t+10}+9y^{16t+11}+9y^{16t+12}+6y^{16t+13}+2y^{16t+14}+y^{16t+15}$\\
$D_{31t+24}$ &
$1+9y^{16t+11}+9y^{16t+12}+6y^{16t+13}+6y^{16t+14}+y^{16t+15}$\\
$D_{31t+25}$ &
$1+7y^{16t+11}+7y^{16t+12}+6y^{16t+13}+6y^{16t+14}+3y^{16t+15}+2y^{16t+16}$\\
$D_{31t+27}$ &
$1+4y^{16t+12}+9y^{16t+13}+9y^{16t+14}+6y^{16t+15}+y^{16t+16}+y^{16t+17}+y^{16t+18}$\\
$D_{31t+28}$ &
$1+9y^{16t+13}+9y^{16t+14}+6y^{16t+15}+5y^{16t+16}+y^{16t+17}+y^{16t+18}$\\
$D_{31t+29}$ &
$1+5y^{16t+13}+5y^{16t+14}+10y^{16t+15}+9y^{16t+16}+y^{16t+17}+y^{16t+18}$\\
$D_{31t+30}$ &
$1+9y^{16t+14}+9y^{16t+15}+5y^{16t+16}+6y^{16t+17}+y^{16t+18}+y^{16t+19}$\\
\hline
$D_{31t}$ &
$1+3y^{16t-2}+9y^{16t-1}+9y^{16t}+6y^{16t+1}+3y^{16t+2}+y^{16t+3}$\\
$D_{31t+16}$ &
$1+6y^{16t+6}+7y^{16t+7}+4y^{16t+8}+6y^{16t+9}+4y^{16t+10}+3y^{16t+11}+y^{16t+12}$\\
\noalign{\hrule height0.8pt}
\end{tabular}
}
\end{center}
\end{table}
%%%%%%%%%%%%%%%%%%%%%%%%%%%%%%%

%%%%%%%%%%%%%%%%%%%%%%%%%%%%%%
\begin{table}[thbp]
\caption{$\det(G(a)G(a)^T)$ for $D_{31t+s}$}
\label{Tab:dim5-det}
\begin{center}
%{\small
{\footnotesize
%{\scriptsize
%{\tiny
\begin{tabular}{c|l}
\noalign{\hrule height0.8pt}
Code & \multicolumn{1}{c}{$\det(G(a)G(a)^T)$}\\
\hline
$D_{31t+3}$ & $196608t^5+61440t^4-1024t^3-1280t^2-48t+1$\\
$D_{31t+4}$ & $196608t^5+69632t^4+7168t^3-32t-1$\\
$D_{31t+5}$ & $196608t^5+102400t^4+20480t^3+1920t^2+80t+1$\\
$D_{31t+7}$ & $196608t^5+192512t^4+66560t^3+9856t^2+640t+15$\\
$D_{31t+11}$&$196608t^5+323584t^4+195584t^3+53888t^2+6704t+301$\\
$D_{31t+19}$&$196608t^5+593920t^4+683008t^3+374656t^2+97904t+9697$\\
$D_{31t+20}$& $196608t^5+684032t^4+944128t^3+647040t^2+220400t+29875$\\
$D_{31t+22}$& $196608t^5+724992t^4+1054720t^3+758400t^2+270000t+38125$\\
$D_{31t+26}$& $196608t^5+856064t^4+1480704t^3+1271168t^2+541296t+91385$\\
\hline
$D_{31t+1}$&
$196608t^5+8192t^4-9216t^3+64t^2+64t-1$ \\
$D_{31t+2}$&
$196608t^5+16384t^4-4096t^3-320t^2+16t+1$ \\
$D_{31t+6}$&
$196608t^5+135168t^4+33792t^3+3840t^2+192t+3$ \\
$D_{31t+8}$&
$196608t^5+212992t^4+82944t^3+14144t^2+1008t+23$ \\
$D_{31t+9}$&
$196608t^5+249856t^4+118784t^3+25728t^2+2448t+81$ \\
$D_{31t+10}$&
$196608t^5+286720t^4+150528t^3+33984t^2+3088t+91$ \\
$D_{31t+12}$&
$196608t^5+356352t^4+240640t^3+73344t^2+9472t+349$ \\
$D_{31t+13}$&
$196608t^5+389120t^4+291840t^3+101568t^2+15792t+803$ \\
$D_{31t+14}$&
$196608t^5+421888t^4+343040t^3+129536t^2+21808t+1167$ \\
$D_{31t+15}$&
$196608t^5+458752t^4+405504t^3+166464t^2+30496t+1797$ \\
$D_{31t+17}$&
$196608t^5+528384t^4+539648t^3+256256t^2+54080t+3575$ \\
$D_{31t+18}$&
$196608t^5+569344t^4+624640t^3+316928t^2+70368t+4555$ \\
$D_{31t+21}$&
$196608t^5+638976t^4+802816t^3+484544t^2+139264t+15015$ \\
$D_{31t+23}$&
$196608t^5+704512t^4+986112t^3+670912t^2+220560t+27811$ \\
$D_{31t+24}$&
$196608t^5+745472t^4+1114112t^3+819264t^2+296032t+41999$ \\
$D_{31t+25}$&
$196608t^5+753664t^4+1132544t^3+832448t^2+298688t+41751$ \\
$D_{31t+27}$&
$196608t^5+819200t^4+1344512t^3+1085376t^2+430352t+66915$ \\
$D_{31t+28}$&
$196608t^5+860160t^4+1495040t^3+1290112t^2+552592t+93971$ \\
$D_{31t+29}$&
$196608t^5+884736t^4+1576960t^3+1390528t^2+606016t+104319$ \\
$D_{31t+30}$&
$196608t^5+913408t^4+1688576t^3+1552448t^2+709760t+129085$ \\
\hline
$D_{31t}$&
$196608t^5-36864t^4-4096t^3+640t^2+16t-1$ \\
$D_{31t+16}$&
$196608t^5+471040t^4+424960t^3+179328t^2+35248t+2589$ \\
\noalign{\hrule height0.8pt}
\end{tabular}
}
\end{center}
\end{table}
%%%%%%%%%%%%%%%%%%%%%%%%%%%%%%%

% %%%%%%%%%%%%%%%%%%%%%%%%%%%%%%%%%%%%%%%%%%%%%%%%%%%%
% \begin{figure}[thb]
% \label{Fig:2}
% \begin{center}
% %{\small
% {\footnotesize
% %{\scriptsize
% %{\tiny
% \[
% \begin{array}{ll}
% M_{19}=
% \left(
% \begin{array}{c}
% 0 0 0 0 0 0 0 1 1 1 1 1 1 1\\
% 0 1 1 1 0 0 1 1 1 0 1 1 1 0\\
% 0 1 0 1 1 1 0 0 1 1 1 1 1 0\\
% 1 0 0 1 1 0 0 1 0 1 0 1 0 1\\
% 1 1 1 0 1 1 0 0 1 0 0 1 0 1
% \end{array}
% \right),
% &
% M_{20}=
% \left(
% \begin{array}{c}
% 000000011111111\\
% 000111101001110\\
% 101011011100101\\
% 110111010001001\\
% 011101100111111
% \end{array}
% \right),
% \\
% M_{22}=
% \left(
% \begin{array}{c}
% 0 0 0 0 0 0 0 0 1 1 1 1 1 1 1 1 1\\
% 1 0 1 1 1 1 0 1 0 1 0 0 1 1 0 1 0\\
% 1 0 1 0 0 1 1 0 0 1 1 1 0 0 0 1 1\\
% 1 1 0 0 1 0 1 1 1 0 1 1 1 0 0 1 0\\
% 1 1 1 1 0 0 1 0 1 1 0 1 1 0 1 0 1
% \end{array}
% \right),
% &
% M_{26}=
% \left(
% \begin{array}{c}
% 0 0 0 0 0 0 0 0 0 0 1 1 1 1 1 1 1 1 1 1 1\\
% 1 0 0 1 1 1 1 1 1 0 1 1 1 1 1 0 0 1 1 0 0\\
% 0 1 1 1 0 1 1 1 1 1 0 1 1 0 1 1 1 1 0 0 0\\
% 0 1 0 0 1 0 0 1 1 1 1 1 0 0 0 1 0 1 0 1 1\\
% 1 1 1 1 0 0 0 0 1 0 1 0 0 0 1 0 1 1 1 1 0
% \end{array}
% \right)
% \end{array}
% \]
% }
% \end{center}
% \caption{$M_{i}$ $(i=19,20,22,26)$}
% \end{figure}
% %%%%%%%%%%%%%%%%%%%%%%%%%%%%%%%%%%%%%%%%%%%%%%%%%%%%

%%%%%%%%%%%%%%%%%%%%%%%%%%%%%%
\begin{table}[thbp]
\caption{Matrices $M_{i}$ $(i=19,20,22,26)$}
\label{Tab:dim4-GM0} 
\begin{center}
%{\small
{\footnotesize
%{\scriptsize
%{\tiny
\begin{tabular}{c|l}
\noalign{\hrule height0.8pt}
$M_{19}$&
00000001111111,
01110011101110,
01011100111110,\\&
10011001010101,
11101100100101
\\
\hline
$M_{20}$&
000000011111111,
000111101001110,
101011011100101,\\&
110111010001001,
011101100111111
\\
\hline
$M_{22}$&
00000000111111111,
10111101010011010,
10100110011100011,\\&
11001011101110010,
11110010110110101
\\
\hline
$M_{26}$&
000000000011111111111,
100111111011111001100,
011101111101101111000,\\&
010010011111000101011,
111100001010001011110
\\
\noalign{\hrule height0.8pt}
\end{tabular}
}
\end{center}
\end{table}
%%%%%%%%%%%%%%%%%%%%%%%%%%%%%%%

Similarly, 
for the parameters
$[31t+s,5,16t+u]$  $(t \ge 1)$, 
where 
\[
(s,u)= (0,-2), (1,-1), (2,-1), (8, 2), (10,3), (12, 4),
\]
and 
$[31t+s,5,16t+u]$ $(t \ge 0)$, 
where 
\begin{align*}
(s,u)=&
( 6,1),
( 9,3),
(13,5),
(14,5),
(15,6),
(16,6),
(17,7),
(18,7),
(21,9),
\\&
(23,10),
(24,11),
(25,11),
(27,12),
(28,13),
(29,13),
(30,14),
\end{align*}
we found the codes $C_{31t+s}=C(a)$.
%($s\in \{0,\ldots,30\} \setminus \{3,4,5,7,11,19,20,22,26\}$).
These codes have  minimum weight one or two less than
the largest possible minimum weight in the bound~\eqref{eq:b5}.
For these codes, we list the vectors $a$, the weight enumerators $W$ and
the determinants $\det(G(a)G(a)^T)$ in 
Tables~\ref{Tab:dim5}, \ref{Tab:dim5-WE} and~\ref{Tab:dim5-det},
respectively.
In Table~\ref{Tab:dim5}, 
we denote $t-1,t+1$ by $t_-,t_+$, respectively.
% From the construction, the dual distances are
% $2$ except $(s,u)=(6,1), (14,5), (19,7), (21,9), (24,11)$ $(t=0)$.
% For the exception, the dual distances are $3$.
Since there is an LCD $[n,5,d]$ code for
$(n,d)=( 8, 2),(10, 3)$ and $(12, 4)$~\cite{bound} and~\cite{HS},
we have the following:

\begin{prop}
If $n \equiv 1,9,13,15,17,21,23,24,25,27,28,29,30 \pmod{31}$ 
and $n \ge 5$, then
\[
d(n,5)=\left\lfloor \frac{16n}{31}\right\rfloor \text{ or }
\left\lfloor \frac{16n}{31}\right\rfloor-1.
\]
If $n \equiv 2,6,8,10,14,18 \pmod{31}$
and $n \ge 5$, then
\[
d(n,5)=\left\lfloor \frac{16n}{31}\right\rfloor-1 \text{ or }
\left\lfloor \frac{16n}{31}\right\rfloor-2.
\]
If $n \equiv 12 \pmod{31}$ and $n \ge 5$, then
\[
d(n,5)=\left\lfloor \frac{16n}{31}\right\rfloor-2 \text{ or }
\left\lfloor \frac{16n}{31}\right\rfloor-3.
\]
If $n \equiv 0,16 \pmod{31}$ and $n \ge 5$, then
\[
d(n,5)=
\left\lfloor \frac{16n}{31}\right\rfloor,
\left\lfloor \frac{16n}{31}\right\rfloor -1
\text{ or }
\left\lfloor \frac{16n}{31}\right\rfloor-2.
\]
\end{prop}

For 
\begin{equation}\label{eq:dim5}
(n,d)=
(25,12),
(27,13),
(28,14),
(29,14),
(30,15),
\end{equation}
in order to verify that there is no LCD $[n,5,d]$
code,
our computer calculation completed a classification of 
(unrestricted) $[n,5,d]$ codes
by the method given in Section~\ref{Sec:k4}.
In Table~\ref{Tab:dim5-small},
we list $N_{n,5,d}$ and 
$N_{n,k,d'}$ ($k=2,3,4, d'=d,d+1,d+2,d+3$) 
for $(n,d)$ in~\eqref{eq:dim5}.
In order to give generator matrices 
$
\left(
\begin{array}{ccccc}
I_5 &  M_{n,i}  \\
\end{array}
\right)
$
of all inequivalent $[n,5,d]$ codes,
we only list the five rows $m_1,m_2,m_3,m_4,m_5$ of $M_{n,i}$ in Table~\ref{Tab:dim5-GM}.
Similar to Table~\ref{Tab:dim4-GM},
the sequences $(m_1,m_2,m_3,m_4,m_5)$
are written in octal using $0=(000),$ $1=(001),\ldots,\ 7=(111)$,
together with $a=(0)$ and $b=(1)$.

%%%%%%%%%%%%%%%%%%%%%%%%%%%%%%
\begin{table}[thbp]
\caption{Numbers of $[n,5,d]$ codes for $(n,d)$ in~\eqref{eq:dim5}}
\label{Tab:dim5-small} 
\begin{center}
%{\small
{\footnotesize
%{\scriptsize
%{\tiny
\begin{tabular}{l|ccccccc}
\noalign{\hrule height0.8pt}
\multicolumn{1}{c|}{$(n,d)$ } 
&$(25,12)$&$(27,13)$&$(28,14)$&$(29,14)$&$(30,15)$\\
\hline
$N_{n,5,d}$     & 8& 1& 1& 9& 1\\
\hline
$N_{n-1,4,d}$   &11& 2& 1&13& 1\\
%$N_{n-1,4,d+1}$ & -& -& -& -& -\\
\hline
$N_{n-2,3,d}$   &16&13& 7&28& 6\\
$N_{n-2,3,d+1}$ & -& 1& -& 1& 1\\
%$N_{n-2,3,d+2}$ &  &  &  & & \\
%$N_{n-2,3,d+3}$ &  &  &  & & \\
\hline
$N_{n-3,2,d}$   &15&15&15&21&15\\
$N_{n-3,2,d+1}$ & 6&10& 6&10&10\\
$N_{n-3,2,d+2}$ & 3& 3& 3& 6& 3\\
$N_{n-3,2,d+3}$ & -& 1& -& 1& 1\\
\noalign{\hrule height0.8pt}
\end{tabular}
}
\end{center}
\end{table}
%%%%%%%%%%%%%%%%%%%%%%%%%%%%%%%

%%%%%%%%%%%%%%%%%%%%%%%%%%%%%%
\begin{table}[thbp]
\caption{Generator matrices of $[n,5,d]$ codes for $(n,d)$ in~\eqref{eq:dim5}}
\label{Tab:dim5-GM} 
\begin{center}
%{\small
{\footnotesize
%{\scriptsize
%{\tiny
\begin{tabular}{c|l}
\noalign{\hrule height0.8pt}
$M_{25,1}$&273153117315434776000007777760547b\\
$M_{25,2}$&546617117315434776000007777760547b\\
$M_{25,3}$&323615531466634773000007777760547b\\
$M_{25,4}$&465630761547437636000007777266633a\\
$M_{25,5}$&466530761547437636000007777266633a\\
$M_{25,6}$&236363174077037770000007776616632b\\
$M_{25,7}$&073663166617037336000007776616632b\\
$M_{25,8}$&263531530747437336000007776616632b\\
$M_{27,1}$&263663176303615761714000037776375746aa\\
$M_{28,1}$&43663307741547434377600000377773721733a\\
$M_{29,1}$&4365154676031760775570000001777732725475\\
$M_{29,2}$&1627171457614660775670000001777732725475\\
$M_{29,3}$&7303615454636660775554000001777732725475\\
$M_{29,4}$&7155415454770660774374000001777732725475\\
$M_{29,5}$&7164314676154360774374000001777732725475\\
$M_{29,6}$&4367714654754360774176000001777732725475\\
$M_{29,7}$&5317606654746630774155400001777732725475\\
$M_{29,8}$&4353606654636630774155400001777732725475\\
$M_{29,9}$&5317606654636630774155400001777732725475\\
$M_{30,1}$&45571433147570361760776000001777747566530bb\\
\noalign{\hrule height0.8pt}
\end{tabular}
}
\end{center}
\end{table}
From the above classification,
our computer calculation shows the following result.

\begin{prop}
There is no LCD $[n,5,d]$ code for
\[
(n,d)=
(25,12),
(27,13),
(28,14),
(29,14),
(30,15).
\]
\end{prop}

%%%%%%%%%%%%%%%%%%%%%%%%%%%%%%%%%%%%%%%%%%%%%%%%%
\section{Largest minimum weights}

The largest minimum weights $d(n,k)$ among all LCD $[n,k]$
codes were determined in~\cite{bound} and~\cite{HS}
for $n \le 12$ and $13 \le n \le 16$, respectively.
In this section, we extend the results to lengths up to $24$.

It is trivial that $d(n,n)=1$.
It is known~\cite{DKOSS} that 
$(d(n,1),d(n,n-1))=(n,2)$ and $(n-1,1)$ 
if  $n$ is odd and even, respectively.
It was shown in~\cite{bound} that
$d(n,2) =
\lfloor \frac{2n}{3}\rfloor$ if  $n \equiv 1,2,3,4 \pmod{6}$, and
$\lfloor \frac{2n}{3} \rfloor-1$ otherwise for $n \ge 2$.
It was shown in~\cite{HS} that
$d(n,3) =\left\lfloor \frac{4n}{7}\right\rfloor$
if $n \equiv 3,5 \pmod {7}$ and 
$\left\lfloor \frac{4n}{7}\right\rfloor-1$ otherwise
for $n \ge 3$.
In addition, by~\cite[Theorem~3]{bound}, 
$d(n,n-2)=2$ for $n \ge 4$, $d(n,n-3)=2$ for $n \ge 8$
and $d(n,n-4)=2$ for $n \ge 16$.
Hence, we only consider the values $d(n,k)$
for $4 \le k \le n-5$ and $17 \le n \le 24$.

We describe how our computer calculation determined
the values $d(n,k)$.
Let $d_{\text{all}}(n,k)$ denote the largest minimum weight
among all (unrestricted) $[n,k]$ codes.
One can find the current information on $d_{\text{all}}(n,k)$ 
in~\cite{G}.
For the following pairs
\begin{equation}\label{eq:P0}
\begin{split}
(n,k)=&
(17,5),
(17,8),
(18,6),
(18,9),
(19,7),
(20,4),
(20,8),
(20,10),
\\&
(21,5),
(21,9),
(22,10),
(22,11),
(23,6),
(23,11),
(23,12),
\\&
(23,14),
(24,5),
(24,7),
(24,12),
(24,14),
\end{split}
\end{equation}
classifications of  $[n,k,d_{\text{all}}(n,k)]$ 
codes are known (see Table~\ref{Tab:known} for the references).
Using the classifications, 
we determined the number $N_L$ of all inequivalent LCD $[n,k,d_{\text{all}}(n,k)]$  codes.
Along with $N_L$,
the number $N$ of all inequivalent $[n,k,d_{\text{all}}(n,k)]$ codes
is listed in Table~\ref{Tab:known}.

%%%%%%%%%%%%%%%%%%%%%%%%%%%%%%
\begin{table}[thbp]
\caption{Known classification of $[n,k,d_{\text{all}}(n,k)]$ codes}
\label{Tab:known}
\begin{center}
%{\small
{\footnotesize
%{\scriptsize
%{\tiny
\begin{tabular}{c|c|cc||c|c|cc}
\noalign{\hrule height0.8pt}
$[n,k,d_{\text{all}}(n,k)]$  & Reference & $N$ &$N_L$ &
$[n,k,d_{\text{all}}(n,k)]$  & Reference & $N$ &$N_L$ \\
\hline
$[17,5,8]$  &\cite{J}       &  5 &  0 &$[22,10,8]$ &\cite{DE}      &  1 & 0 \\
$[17,8,6]$  &\cite{Simonis} &  1 &  1 &$[22,11,7]$ &\cite{BH}      &  1 & 0 \\
$[18,6,8]$  &\cite{DE}      &  2 &  0 &$[23,6,10]$ &\cite{J}       & 29 &14 \\
$[18,9,6]$  &\cite{Simonis} &  1 &  0 &$[23,11,8]$ &\cite{DE}      &  1 & 0 \\
$[19,7,8]$  &\cite{DE}      &  1 &  0 &$[23,12,7]$ & (see~\cite{P})&  1 & 0 \\
$[20,4,10]$ &\cite{J}       &  3 &  1 &$[23,14,5]$ &\cite{Simonis2}&  1 & 0 \\
$[20,8,8]$  &\cite{DE}      &  1 &  0 &$[24,5,12]$ &\cite{J}       &  1 & 0 \\
$[20,10,6]$ &\cite{GO}      &1682&601 &$[24,7,10]$ &\cite{J}       &  6 & 0 \\
$[21,5,10]$ &\cite{Tilborg} &  2 &  0 &$[24,12,8]$ &\cite{Snover}  &  1 & 0 \\
$[21,9,8]$  &\cite{DE}      &  1 &  0 &$[24,14,6]$ &\cite{J}       &  1 & 0 \\
\noalign{\hrule height0.8pt}
\end{tabular}
}
\end{center}
\end{table}
%%%%%%%%%%%%%%%%%%%%%%%%%%%%%%%

For the following pairs
\begin{equation}\label{eq:P1}
\begin{split}
(n,k)=&
(17,  6),(17, 11),(18,  5),(18,  7),(19,  8),(19, 13),(20,  7),(20,  9),
\\&
(21,  8),(21, 10),(21, 11),(21, 15),(22,  4),(22,  9),(23,  4),(23,  5),
\\&
(23,  9),(23, 10),(23, 13),(23, 17),(24, 11),
\end{split}
\end{equation}
we completed the classifications of  $[n,k,d_{\text{all}}(n,k)]$ codes
by the method given in Section~\ref{Sec:k4}.
The number $N$ of all inequivalent $[n,k,d_{\text{all}}(n,k)]$ codes
is listed in Table~\ref{Tab:New}.
From the classification, we know that 
there is no LCD $[n,k,d_{\text{all}}(n,k)]$ code for
$(n,k)$ listed in~\eqref{eq:P1}.
For each of the parameters, 
all inequivalent codes can be obtained electronically from
\url{http://yuki.cs.inf.shizuoka.ac.jp/lcd2/}.
In addition, for the following pairs
\begin{equation}\label{eq:P2}
(n,k)=(20,8), (21,9), (22,10), (23,11), (24,12),
\end{equation}
we also 
completed the classifications of  $[n,k,d_{\text{all}}(n,k)-1]$ codes
by the method given in Section~\ref{Sec:k4}.
The number $N'$ of all inequivalent $[n,k,d_{\text{all}}(n,k)-1]$ codes
is listed in Table~\ref{Tab:New-2}.
From the classification, we know that
there is no LCD $[n,k,d_{\text{all}}(n,k)-1]$ code
for $(n,k)$ listed in~\eqref{eq:P2}.
For each of the parameters, 
all inequivalent codes can be obtained electronically from
\url{http://yuki.cs.inf.shizuoka.ac.jp/lcd2/}.

%%%%%%%%%%%%%%%%%%%%%%%%%%%%%%
\begin{table}[thbp]
\caption{Classification of $[n,k,d_{\text{all}}(n,k)]$ codes}
\label{Tab:New}
\begin{center}
%{\small
{\footnotesize
%{\scriptsize
%{\tiny
\begin{tabular}{c|c||c|c||c|c}
\noalign{\hrule height0.8pt}
$[n,k,d_{\text{all}}(n,k)]$  & $N$ &
$[n,k,d_{\text{all}}(n,k)]$  & $N$ &
$[n,k,d_{\text{all}}(n,k)]$  & $N$ \\
\hline
$[17,  6, 7]$&    3 &$[20,  9, 7]$&    1 &$[23,  4,12]$&    1 \\
$[17, 11, 4]$&   40 &$[21,  8, 8]$&   13 &$[23,  5,11]$&    1 \\
$[18,  5, 8]$&   39 &$[21, 10, 7]$&    1 &$[23,  9, 8]$&40289 \\
$[18,  7, 7]$&    2 &$[21, 11, 6]$&  739 &$[23, 10, 8]$&    9 \\
$[19,  8, 7]$&    1 &$[21, 15, 4]$&   16 &$[23, 13, 6]$&    8 \\
$[19, 13, 4]$&   25 &$[22,  4,11]$&    2 &$[23, 17, 4]$&    9 \\
$[20,  7, 8]$&   27 &$[22,  9, 8]$&   10 &$[24, 11, 8]$&   10 \\
\noalign{\hrule height0.8pt}
\end{tabular}
}
\end{center}
\end{table}
%%%%%%%%%%%%%%%%%%%%%%%%%%%%%%%

%%%%%%%%%%%%%%%%%%%%%%%%%%%%%%
\begin{table}[thbp]
\caption{Classification of $[n,k,d_{\text{all}}(n,k)-1]$ codes}
\label{Tab:New-2}
\begin{center}
%{\small
{\footnotesize
%{\scriptsize
%{\tiny
\begin{tabular}{c|c||c|c}
\noalign{\hrule height0.8pt}
$[n,k,d_{\text{all}}(n,k)-1]$  & $N'$ &
$[n,k,d_{\text{all}}(n,k)-1]$  & $N'$ \\
\hline
$[20, 8,7]$ & 34 &$[23,11,7]$ & 16 \\
$[21, 9,7]$ & 21 &$[24,12,7]$ & 12 \\
$[22,10,7]$ & 16 & & \\
\noalign{\hrule height0.8pt}
\end{tabular}
}
\end{center}
\end{table}
%%%%%%%%%%%%%%%%%%%%%%%%%%%%%%%

Finally, 
for the pairs $(n,k)$ listed in~\eqref{eq:P0}, \eqref{eq:P1} 
and~\eqref{eq:P2}
except $(17,8)$, $(20,4)$, $(20,10)$ and $(23,6)$,
we found an LCD $[n,k,d]$ code, where $d$ is the minimum
weight listed in Table~\ref{Tab:maxd}.
Then we determined the values $d(n,k)$
for $4 \le k \le n-5$ and $17 \le n \le 24$.
For each of the parameters listed in the table, 
an LCD code can be obtained electronically from
\url{http://yuki.cs.inf.shizuoka.ac.jp/lcd2/}.

%%%%%%%%%%%%%%%%%%%%%%%%%%%%%%
\begin{table}[thbp]
\caption{$d(n,k)$ for $17 \le n \le 24$}
\label{Tab:maxd}
\begin{center}
%{\small
{\footnotesize
%{\scriptsize
%{\tiny
\begin{tabular}{c|ccccccccccccccccccccc}
\noalign{\hrule height0.8pt}
$n\backslash k$  & 4& 5& 6& 7& 8& 9&10&11&12&13&14&15&16&17&18&19 \\
\hline
17& 8& 7& 6& 6& 6& 5  & 4& 3& 3&  &  &    \\
18& 8& 7& 7& 6& 6& 5  & 4& 4& 4& 3&  &    \\
19& 9& 8& 8& 7& 6& 6  & 5& 4& 4& 3& 3&    \\
20&10& 9& 8& 7& 6& 6  & 6& 5& 4& 4& 4& 3  \\
21&10& 9& 8& 8& 7& 6  & 6& 5& 5& 4& 4& 3& 3  \\
22&10&10& 9& 8& 8& 7  & 6& 6& 6& 5& 4& 4& 4 &3   \\
23&11&10& 9& 9& 8& 7  & 7& 6& 6& 5& 4& 4& 4 &3 &3   \\
24&12&11&10& 9& 8& 8  & 8& 7& 6& 6& 5& 4& 4 &4 &4 & 3 \\
\noalign{\hrule height0.8pt}
\end{tabular}
}
\end{center}
\end{table}
%%%%%%%%%%%%%%%%%%%%%%%%%%%%%%%

%%%%%%%%%%%%%%%%%%%%
\bigskip
\noindent
{\bf Acknowledgment.}
This work was supported by JSPS KAKENHI Grant Number 15H03633.
In this work, the supercomputer of ACCMS, Kyoto University
was partially used.

%%%%%%%%%%%%%%%%%%%  References  %%%%%%%%%%%%%%%%%%%%%%%%

%%%%%%%%%%%%%%%%%%%%%%%%%%%%%%%%%

\begin{landscape}
%%%%%%%%%%%%%%%%%%%%%%%%%%%%%%
\begin{table}[thb]
\caption{LCD codes of dimension 5}
\label{Tab:dim5}
\begin{center}
%{\small
{\footnotesize
%{\scriptsize
%{\tiny
\begin{tabular}{c|l}
\noalign{\hrule height0.8pt}
Code & \multicolumn{1}{c}{$a$}\\
\hline
$D_{31t+3}$ &
%$(t,t,t,t,t,t,t,t,t,t,t,t,t,t,t,t,t,t,t,t,t,t,t-1,t,t,t+1,t+1,t-1,t-1,t-1,t)$\\
$(t,t,t,t,t,t,t,t,t,t,t,t,t,t,t,t,t,t,t,t,t,t,t_{-},t,t,t_{+},t_{+},t_{-},t_{-},t_{-},t)$\\
$D_{31t+4}$ &
%$(t,t,t,t,t,t,t,t,t,t,t,t,t,t,t,t,t,t,t,t,t,t,t,t,t,t,t,t,t-1,t,t)$\\
$(t,t,t,t,t,t,t,t,t,t,t,t,t,t,t,t,t,t,t,t,t,t,t,t,t,t,t,t,t_{-},t,t)$\\
$D_{31t+5}$ &
$(t,t,t,t,t,t,t,t,t,t,t,t,t,t,t,t,t,t,t,t,t,t,t,t,t,t,t,t,t,t,t)$\\
$D_{31t+7}$ &
%<t,t,t,t,t,t,t,t,t,t,t,t,t,t,t+1,t,t,t,t,t,t,t+1,t,t,t,t+1,t+1,t,t,t-1,t-1>
$(t,t,t,t,t,t,t,t,t,t,t,t,t,t,t_+,t,t,t,t,t,t,t_+,t,t,t,t_+,t_+,t,t,t_-,t_-)$ \\
$D_{31t+11}$&
%<t,t,t,t,t,t,t,t,t,t,t,t+1,t+1,t+1,t+1,t-1,t,t,t+
%1,t+1,t+1,t,t,t,t,t,t,t,t,t,t>
$(t,t,t,t,t,t,t,t,t,t,t,t_+,t_+,t_+,t_+,t_-,t,t,t_+,t_+,t_+,t,t,t,t,t,
t,t,t,t,t)$\\
$D_{31t+19}$ &
%<t+1,t+1,t+1,t+1,t+1,t+1,t+1,t+1,t+1,t+1,t+1,t
%+1,t+1,t+1,t+1,t+1,t,t,t,t,t+1,t,t,t,t+1,t,t-
%1,t-1,t-1,t,t-1>
$(t_+,t_+,t_+,t_+,t_+,t_+,t_+,t_+,t_+,t_+,t_+,t_+,t_+,t_+,t_+,t_+,t,t,t,
t,t_+,t,t,t,t_+,t,t_-,t_-,t_-,t,t_-)$ \\
$D_{31t+20}$&
%t+1,t+1,t+1,t+1,t+1,t+1,t+1,t+1,t+1,t+1,t+1,t
%+1,t+1,t+1,t+1,t-1,t+1,t+1,t+1,t,t+1,t,t,t-1,
%t+1,t,t,t-1,t,t-1,t-1
$(t_+,t_+,t_+,t_+,t_+,t_+,t_+,t_+,t_+,t_+,t_+,t_+,t_+,t_+,t_+,t_-,t_+,t_+,t_+,t,t_+,t,t,t_-,t_+,t,t,t_-,t,t_-,t_-)$\\
$D_{31t+22}$&
%$(t+1,t+1,t+1,t+1,t+1,t+1,t+1,t+1,t+1,t+1,t+1,t+1,t+1,t+1,t+1,t+1,t+1,t+1,t+1,t,t+1,t,t,t-1,t+1,t,t,t-1,t,t-1,t-1)$\\
$(t_{+},t_{+},t_{+},t_{+},t_{+},t_{+},t_{+},t_{+},t_{+},t_{+},t_{+},t_{+},t_{+},t_{+},t_{+},t_{+},t_{+},t_{+},t_{+},t,t_{+},t,t,t_{-},t_{+},t,t,t_{-},t,t_{-},t_{-})$\\
$D_{31t+26}$&
%$(t+1,t+1,t+1,t+1,t+1,t+1,t+1,t+1,t+1,t+1,t+1,t+1,t+1,t+1,t+1,t+1,t+1,t+1,t+1,t+1,t+1,t+1,t,t-1,t+1,t,t+1,t-1,t+1,t-1,t-1)$\\
$(t_{+},t_{+},t_{+},t_{+},t_{+},t_{+},t_{+},t_{+},t_{+},t_{+},t_{+},t_{+},t_{+},t_{+},t_{+},t_{+},t_{+},t_{+},t_{+},t_{+},t_{+},t_{+},t,t_{-},t_{+},t,t_{+},t_{-},t_{+},t_{-},t_{-})$\\
\hline
$D_{31t+1}$ &
$(t,t,t,t,t,t,t,t,t,t,t,t,t,t,t_-,t,t,t,t,t,t,t,t_-,t_-,t,t_+,t_+,t_-,t_-,t_-,t)$ \\
$D_{31t+2}$ &
$(t,t,t,t,t,t,t,t,t,t,t,t,t,t,t,t,t,t,t,t,t,t,t,t,t,t,t_-,t,t_-,t_-,t)$ \\
$D_{31t+6}$ &
$(t,t,t,t,t,t,t,t,t,t,t,t,t,t,t,t,t,t,t,t,t,t,t,t,t,t,t,t,t_+,t,t)$ \\
$D_{31t+8}$ &
$(t,t,t,t,t,t,t,t_+,t,t_+,t,t,t+2,t_-,t,t,t,t,t,t,t,t,t,t,t,t,t,t,t,t,t)$ \\
$D_{31t+9}$ &
$(t,t,t,t,t,t_+,t_+,t,t,t_+,t_+,t,t,t,t,t,t,t,t,t,t,t,t,t,t,t,t,t,t,t,t)$ \\
$D_{31t+10}$ &
$(t,t,t,t,t,t,t_+,t_+,t,t+2,t,t,t+2,t_-,t,t,t,t,t,t,t,t,t,t,t,t,t,t,t,t,t)$ \\
$D_{31t+12}$ &
$(t,t,t,t,t,t_+,t_+,t_+,t,t+2,t,t_+,t+2,t_-,t,t,t,t,t,t,t,t,t,t,t,t,t,t,t,t,t)$ \\
$D_{31t+13}$ &
$(t,t,t,t_+,t,t_+,t+2,t,t_+,t_+,t,t_+,t_+,t,t,t,t,t,t,t,t,t,t,t,t,t,t,t,t,t,t)$ \\
$D_{31t+14}$ &
$(t,t,t,t_+,t,t_+,t_+,t_+,t_+,t_+,t_+,t,t_+,t_+,t,t,t,t,t,t,t,t,t,t,t,t,t,t,t,t,t)$ \\
$D_{31t+15}$ &
$(t,t,t,t+2,t_+,t_+,t_+,t,t_+,t_+,t_+,t_+,t_+,t,t,t,t,t,t,t,t,t,t,t,t,t,t,t,t,t,t)$ \\
$D_{31t+17}$ &
$(t,t,t,t+2,t_+,t_+,t_+,t_+,t_+,t_+,t+2,t,t_+,t_+,t,t,t,t,t,t,t,t,t,t,t,t,t,t,t,t,t)$ \\
$D_{31t+18}$ &
$(t,t,t,t+2,t,t+2,t+2,t,t_+,t_+,t_+,t_+,t+2,t_+,t,t,t,t,t,t,t,t,t,t,t,t,t,t,t,t,t)$ \\
$D_{31t+21}$ &
$(t_+,t_+,t_+,t_+,t_+,t_+,t_+,t_+,t_+,t_+,t_+,t_+,t_+,t_+,t,t,t,t,t,t,t,t,t_+,t,t,t_+,t,t,t,t,t)$ \\
$D_{31t+23}$ &
$(t_+,t_+,t_+,t_+,t_+,t_+,t_+,t_+,t_+,t_+,t_+,t_+,t_+,t_+,t,t,t,t,t,t,t,t_+,t_+,t,t,t_+,t_+,t,t,t,t)$ \\
$D_{31t+24}$ &
$(t_+,t_+,t_+,t_+,t_+,t_+,t_+,t_+,t_+,t_+,t_+,t_+,t_+,t_+,t,t,t,t,t,t_+,t_+,t,t_+,t,t_+,t_+,t,t,t,t,t)$ \\
$D_{31t+25}$ &
$(t_+,t_+,t_+,t_+,t_+,t_+,t_+,t_+,t_+,t_+,t_+,t_+,t_+,t_+,t,t,t,t,t,t,t_+,t,t_+,t_+,t_+,t_+,t,t_+,t,t,t)$ \\
$D_{31t+27}$ &
$(t_+,t_+,t_+,t_+,t_+,t_+,t_+,t_+,t_+,t_+,t_+,t_+,t_+,t_+,t,t,t,t,t,t,t_+,t_+,t_+,t_+,t_+,t_+,t_+,t_+,t,t,t)$ \\
$D_{31t+28}$ &
$(t_+,t_+,t_+,t_+,t_+,t_+,t_+,t_+,t_+,t_+,t_+,t_+,t_+,t_+,t,t,t,t_+,t_+,t,t_+,t,t_+,t_+,t_+,t_+,t,t_+,t_+,t,t)$ \\
$D_{31t+29}$ &
$(t_+,t_+,t_+,t_+,t_+,t_+,t_+,t_+,t_+,t_+,t_+,t_+,t_+,t_+,t,t,t,t,t,t_+,t_+,t_+,t_+,t_+,t_+,t_+,t_+,t_+,t_+,t,t)$ \\
$D_{31t+30}$ &
$(t_+,t_+,t_+,t_+,t_+,t_+,t_+,t_+,t_+,t_+,t_+,t_+,t_+,t_+,t,t_+,t,t_+,t_+,t,t_+,t_+,t_+,t,t_+,t_+,t_+,t_+,t,t,t_+)$ \\
\hline
$D_{31t}$ &
$(t,t,t,t,t,t,t,t,t,t,t,t,t,t,t,t,t,t,t,t,t,t,t,t,t,t_-,t_-,t,t_-,t_-,t_-)$ \\
$D_{31t+16}$ &
$(t,t,t,t,t,t,t,t,t,t,t,t_+,t,t+2,t+2,t,t_+,t,t_+,t_+,t_+,t_+,t,t,t,t,t_+,t,t,t,t)$ \\
\noalign{\hrule height0.8pt}
\end{tabular}
}
\end{center}
\end{table}
%%%%%%%%%%%%%%%%%%%%%%%%%%%%%%%

\end{landscape}


\begin{thebibliography}{99}


\bibitem{BH}K. Betsumiya and M. Harada,
Binary optimal odd formally self-dual codes,
{\sl Des.\ Codes Cryptogr.}
{\bf 23} (2001), 11--21.

\bibitem{Magma}W. Bosma, J. Cannon and C. Playoust,
The Magma algebra system I: The user language,
{\sl J. Symbolic Comput.}
{\bf 24} (1997), 235--265.

% \bibitem{CG}
% C. Carlet and S. Guilley,
% Complementary dual codes for counter-measures to side-channel attacks,
% In: E.R. Pinto et al.\ (eds.), 
% Coding Theory and Applications, CIM Series in
% Mathematical Sciences, vol.\ 3, pp.~97--105, Springer, 2014.

% \bibitem{CMTQ}
% C. Carlet, S. Mesnager, C. Tang and Y. Qi,
% Euclidean and Hermitian LCD MDS codes,
% {\sl Des.\ Codes Cryptogr.},
% (to appear),
% arXiv:1702.08033.

\bibitem{mf}
C. Carlet, S. Mesnager, C. Tang and Y. Qi,
New characterization and parametrization of LCD codes,
{\sl IEEE\ Trans.\ Inform.\ Theory},
(to appear),
arXiv:1709.03217.

\bibitem{CMTQ2}
C. Carlet, S. Mesnager, C. Tang, Y. Qi and R. Pellikaan,
Linear codes over $\FF_q$ are equivalent to LCD codes for $q >3$,
{\sl IEEE\ Trans.\ Inform.\ Theory}
{\bf 64}  (2018),  3010--3017.

\bibitem{DE}
S.M. Dodunekov and S.B. Encheva, 
Uniqueness of some linear subcodes of the binary extended Golay
code,
{\sl Problems Inform.\ Transmission}
{\bf 29}  (1993),  38--43.

\bibitem{DKOSS}
S.T. Dougherty, J.-L. Kim, B. Ozkaya, L. Sok and P. Sol\'e,
The combinatorics of LCD codes: linear programming bound and orthogonal
matrices,
{\sl Int.\ J. Inf.\ Coding Theory}
{\bf 4}  (2017),  116--128.

\bibitem{bound}
L. Galvez, J.-L. Kim, N. Lee, Y.G. Roe and B.-S. Won,
Some bounds on binary LCD codes,
{\sl Cryptogr.\ Commun.}
{\bf 10} (2018), 719--728.

\bibitem{G} M. Grassl,
Code tables: Bounds on the parameters of various types of codes,
Available online at 
\url{http://www.codetables.de/},
Accessed on 2018-04-06.

\bibitem{Gr} J.H. Griesmer, 
A bound for error-correcting codes. 
{\sl IBM J. Res.\ Develop.}
{\bf 4} (1960), 532--542. 

\bibitem{GO}T.A. Gulliver and P.R.J. \"Osterg\aa rd,
Binary optimal linear rate $1/2$ codes,
{\sl Discrete Math.}
{\bf 283} (2004), 255--261.


% \bibitem{GOS}C. \"Guneri, B. \"Ozkaya and P. Sol\'e, 
% Quasi-cyclic complementary dual codes,
% {\sl Finite Fields Appl.}
% {\bf 42}  (2016), 67--80.

\bibitem{HS} M. Harada and K. Saito,
Binary linear complementary dual codes,
{\sl Cryptogr.\ Commun.}, (to appear),
arXiv:1802.06985.

\bibitem{J}D.B. Jaffe, 
Optimal binary linear codes of length $\le 30$,
{\sl Discrete Math.}
{\bf 223} (2000),  135--155.

% \bibitem{Jin}L. Jin,
% Construction of MDS codes with complementary duals,
% {\sl IEEE Trans.\ Inform.\ Theory}
% {\bf 63} (2017), 2843--2847.

\bibitem{Massey}J.L. Massey, 
Linear codes with complementary duals,
{\sl Discrete Math.}
{\bf 106/107} (1992), 337--342.

\bibitem{nauty}
B.D. McKay and A. Piperno, 
Practical graph isomorphism, II, 
{\sl J. Symbolic Comput.}
{\bf 60} (2014), 94--112.

% \bibitem{LN} E.R. Lina, Jr.\ and E.G. Nocon,
% On the construction of some LCD codes over finite fields,
% {\sl Manila J. Science}
% {\bf 9} (2016), 67--82.

\bibitem{P}V. Pless, 
Introduction to the theory of error-correcting codes (Third edition), 
John Wiley \& Sons, Inc., New York, 1998.

% \bibitem{S}M. Sendrier,
% Linear codes with complementary duals meet the Gilbert--Varshamov bound,
% {\sl Discrete Math.}
% {\bf 285} (2004), 345--347.

\bibitem{NTL}V. Shoup,
NTL: A Library for doing Number Theory,
Available online at 
\url{http://www.shoup.net/ntl/}.

\bibitem{Simonis} J. Simonis,
The $[18, 9, 6]$ code is unique, 
{\sl Discrete Math.}
{\bf 106/107} (1992) 439--448.

\bibitem{Simonis2} J. Simonis,
The $[23,14,5]$ Wagner code is unique,
{\sl Discrete Math.}
{\bf 213} (2000),  269--282. 

\bibitem{Snover}S.L. Snover, 
The uniqueness of the Nordstrom--Robinson and 
the Golay binary codes, Ph.D.~Thesis,
Michigan State Univ., 1973.

\bibitem{Tilborg}H. van Tilborg, 
On the uniqueness resp.\ nonexistence of certain codes meeting the
Griesmer bound,
{\sl Inform.\ Control}
{\bf 44}  (1980), 16--35. 


% \bibitem{OEIS}N. Sloane and  S. Plouffe,
% The Encyclopedia of Integer Sequences,
% Academic Press, San Diego, CA, 1995
% (Available online at \url{https://oeis.org}).

\end{thebibliography}
\end{document}